\documentclass[11pt]{amsart}
\usepackage{latexsym}
\usepackage{amsthm, amsmath, amssymb,latexsym,xypic, color}

\newtheorem{theorem}{Theorem}[section]

\newtheorem{proposition}[theorem]{Proposition}
\newtheorem{remark}[theorem]{Remark}
\newtheorem{definition}[theorem]{Definition}

\newcommand{\nc}{\newcommand}
\nc{\cC}{{\mathcal C}} \nc{\cO}{{\mathcal O}}
\nc{\cS}{{\mathcal A}} \nc{\cE}{{\mathcal E}}
\nc{\cZ}{{\mathcal Z}}
\nc{\cM}{{\mathcal M}}
\nc{\mR}{{\mathcal R}}
\nc{\mH}{{\mathcal H}}
\nc{\cT}{{\mathcal T}}
\nc{\uR}{\underline{\mathbb R}} \nc{\uC}{\underline{\mathbb C}}
\nc{\bH}{{\mathbb H}} \nc{\bA}{{\mathbb A}} \nc{\bG}{{\mathbb G}}
\nc{\bC}{{\mathbb C}} \nc{\bO}{{\mathbb O}} \nc{\bI}{{\mathbb I}}
\nc{\bB}{{\mathbb B}} \nc{\bY}{{\mathbb Y}} \nc{\bK}{{\mathbb K}}
\nc{\bX}{{\mathbb X}} \nc{\bS}{{\mathbb S}} \nc{\bE}{{\mathbb E}}
\nc{\bF}{{\mathbb F}} \nc{\bZ}{{\mathbb Z}} \nc{\bQ}{{\mathbb Q}}
\nc{\bN}{{\mathbb N}} \nc{\bP}{{\mathbb P}} \nc{\bL}{{\mathbb L}}
\nc{\bM}{{\mathbb M}} \nc{\bT}{{\mathbb T}} \nc{\bW}{{\mathbb W}}
\nc{\bU}{{\mathbb U}} \nc{\bD}{{\mathbb D}} \nc{\bJ}{{\mathbb J}}
\nc{\bV}{{\mathbb V}} \nc{\bbZ}{{\mathbb Z}} \nc{\bR}{{\mathbb R}}
\nc{\co}{{\nabla}} \nc{\cu}{{\overline{\nabla}}}
\nc{\fr}{{\rightarrow}}

\begin{document}
                                %
                                %
                                %
%%%%%%%%%%%%%%%%%%%%%%%%%%%%%%%%%
%% front page
%%%%%%%%%%%%%%%%%%%%%%%%%%%%%%%%%
\title{On the dimension of some real, bounded rank, matrix spaces} %
                                %
                                %
%\author{Andrea Causin  \footnote{Partially supported by
%1) PRIN 2005 {\em ``Spazi di moduli e teorie di Lie"}; 2) Indam
%(GNSAGA); 3) Far 2006 (PV):{\em ``Variet\`{a} algebriche, calcolo
%algebrico, grafi orientati e topologici"}.}}
% \vspace{-1.5mm}
% {\normalsize Dipartimento di Matematica, Universit\`a di Pavia,}\\
% {\normalsize via Ferrata 1, 27100 Pavia, Italy}\\
% {\footnotesize e-mail: {\tt pirola@dimat.unipv.it}}
%\date{}
%qui si introduce la data

%\thanks{}  %qui ci vanno gli acknowledgments
                                %
                                %
\maketitle
                                %
                                %
%                                 %

%\vskip -4mm
                                %
                                %

%%%%%%%%%%%%%%%%%%%%%%%

                                %
                                %
                                %
%%%%%%%%%%%%%%%%%%%%%%%%%%%%%%%%%
%% front page
%%%%%%%%%%%%%%%%%%%%%%%%%%%%%%%%%
                                %
                                %
%%%%%%%%%%%%%%%%%%%%%%%%%%%%%%%%%
%% text
%%%%%%%%%%%%%%%%%%%%%%%%%%%%%%%%%
 \begin{abstract} {\em
\noindent Given $n\in\bN$, let $X$ be either the set of hermitian or real  $n\times n$ matrices of rank at least $n-1$. If $n$ is even, we give a sharp estimate on the maximal dimension of a real vector space $V\subset X\cup\{0\}$. The rusults are obtained, via $K-$theory, by studying a bundle map induced by the adjugation of matrices.} \vskip3mm \vskip 1mm
 {\setlength{\baselineskip}{0.8\baselineskip}

   \par}

\end{abstract}

%%%%%%%%%%%%%%%%%%%%%%%%%%
%%%%%%%%%%%%%%%%%%%%%%%%%%%%%%%%%

%\section*{Introduction}

%[E QUI CHE ACCIDENTI CI METTO??? LO LEGO COL PARAGRAFO SUCCESSIVO...]

\section{Introduction: real subspaces of matrices}

Let $n>0$ be a positive integer and denote by $\cM_n(\bC)$ the space of $n\times n$ complex matrices. Let $H$ and $R$ be respectively the real subspaces of $\cM_n(\bC)$ of hermitian and real matrices:
$$H=\{A\in\cM_n(\bC)|\, \overline{A}= {^t\!A} \},\quad \quad R=\{A\in\cM_n(\bC)|\, \overline{A}=A\},$$ 
$\dim H=\dim R= n^2.$ 
\begin{definition}
For any vector subspace $V\in \cM_n(\bC)$, the {\em minimal rank} of $V$ is the positive integer 
$$m_V=\min_{A\in \cM_n(\bC)\setminus\{0\}} \mbox{ rank }A.$$
\end{definition}
Let now $\mH_m$ and $\mR_m$ be respectively the sets of linear subspaces of $H$ and $R$ having minimal rank $m$.
\begin{definition}
Fixed integers $n>0$ and $0<m\leq n$, set
$$h_{n,m}=\max_{V\in\mH_m}\dim V \quad \mbox{ and } \quad r_{n,m}=\max_{V\in\mR_m}\dim V$$
\end{definition}
Note that $h_{n,n}\leq h_{n,n-1}\leq\dots\leq h_{n,1}=n^2$ and similarly for $r_{n,m}.$
\medskip

Consider the following:
\medskip

\noindent {\bf Problem} {\em Compute or give an estimate of $h_{n,m}$ and $r_{n,m}$.}
\bigskip

For $n\in\bN$, factorise $n=2^{a+4b}(2k+1)$ with $a,b,k\in\bN$, $0\leq a\leq 3$ and define the real and complex Radon-Hurwitz numbers as $\rho(n)=2^a+8b$, $\rho_\bC(n)=2(a+4b)+2$. In \cite{adams1} and \cite{adams2} Adams, Lax and Phillips show, by determining the maximal number of everywhere independent vector fields on a sphere, that $r_{n,n}=\rho(n)$ and $h_{n,n}=\rho_\bC(n/2)+1$. In this paper we prove:
\medskip

\noindent {\bf Theorem} {\em If $n>0$ is an even integer, then 
$$h_{n,n-1}\leq\rho_\bC(n)\quad \mbox{ and }\quad r_{n,n-1}\leq \rho_\bC(n)$$
and, if $8$ divides $n$, then $r_{n,n-1}= \rho_\bC(n)=\rho(n).$}
\medskip

Our interest in dealing with spaces of hermitian matrices of rank bounded from below primarly arised by studying the kernel of the cup-product map $\Lambda^2H^1(X,\bC)\to H^2(X,\bC)$ of a compact K\"ahler variety $X$ without Albanese fibrations (see the author and Pirola in \cite{avvoltoio}).

The problem posed in this paper is, however, a particular formulation of the following more general question: given a set $X$ of (real or complex) matrices verifying some algebraic conditions such as {\em e.g.} symmetry or boundedness of rank, what is the maximal dimension of a linear subspace $V\subset X\cup\{0\}$? 

This question has an interest on its own since it naturally arises in many different settings. In its original formulation (concerning real invertible matrices) it is equivalent to the problem of finding everywhere indepedent vector fields on a sphere and has been solved, as said before, by Adams and others in \cite{adams1,adams2}.
For symmetric invertible matrices it is related to spectral problems and PDE's (see Friedland, Robbin and Sylvester \cite{crossing}). In the case of real rectangular matrices of maximal rank it is related to the
existence of bilinear nonsingular maps $\bR^k\times \bR^n\to\bR^m$ generalizing the multiplication map of the classical division algebras over $\bR$ (see Hurwitz \cite{hurwitz}, Berger and Friedland \cite{berger}, Lam and Yiu \cite{lamyiu}). In the case of constant rank matrices it is related to the geometric dimension of vector bundles over the projective space, hence to immersion poblems (see Adams \cite{adams3}, Beasley \cite{beasley}, the autor and Pirola \cite{fagiano}, Landsberg and Ilic \cite{landsberg}, Meshulam \cite{meshulam}, Rees \cite{rees1}, Westwick \cite{westwick}). 
\medskip

The proof of our theorem is based on the following considerations: if $A\in\cM_n(\bC)$ is a real or hermitian matrix of rank $\geq n-1$, then its adjugate matrix $A^\ast$ is not zero and $\psi(A)=A+i{^t\, \!\overline{A^\ast}}$ is invertible. When $n$ is even, $\psi(-A)=-\psi(A)$ and, if $V$ is any real linear space of such matrices, the map $\psi:V\setminus\{0\}\to GL_n(\bC)$ induces an isomormphism between the trivial vector bundle $\uC^n$ and $n$ times the complex tautological bundle over $\bP(V)$. This isomorphism, when it is read in the ring of reduced $K-$theory of $\bP(V)$, gives an algebraic relation between $n$ and $\dim V$ which implies the stated estimate. 

The paper is therefore organized in two sections: in the first one we recall the principal results about the structure of the $K-$theory ring for the projective space and we give a method for relating linear spaces of matrices to bundle maps. In the second part we prove our statements about adjugate matrices and conclude the proof of the theorem.

\subsection*{Acknowledgements} The author would like to express his gratitude to S. Friedland for the many fruitful discussions on these subjects. The author is also deeply grateful to G. P. Pirola for the uncountable number of teachings he gave and for his constant interest in the author's work.

\section{Preliminary statements}
For the main cited results about K-theory refer to \cite{adams1} or to \cite{husemoller}; for the definition and properties of odd maps in this setting, the main references are \cite{crossing} and \cite{pacific}.

\subsection{K-theory and Radon-Huwitz numbers}

Let $\bR\bP^{d-1}$ be the real projective space of dimension $d-1$, $\xi$ its tautological vector bundle and $\underline{\bR}^k$, $\underline{\bC}^k$ the real and complex trivial bundles of rank $k$. Denote by $\xi_\bC=\xi\otimes\bC$ the complexification of $\xi$ and remark that 
$\xi_\bC=\xi\oplus\xi$ and $\xi\otimes\xi=\underline{\bR}.$

Recall that the ring of reduced complex K-theory $\widetilde{K}_\bC(\bR\bP^{d-1})$ is the ring of formal differences $[E]-[\bC^k]$ of (isomorphism classes of) complex vector bundles over $\bR\bP^{d-1}$ such that $E$ has rank $k$.

\begin{proposition}\label{anello} The ring $\widetilde{K}_\bC(\bR\bP^{d-1})$ is isomorphic to the polynomial ring $\bZ[\mu]$ with the following relations 
\begin{eqnarray}
\label{rel1}\mu^2=-2\mu & \\ 
\label{rel2} \mu^{g(d)+1}=0, & g(d)=\mbox{integer part of }\frac{d-1}{2}.
\end{eqnarray}
The isomorphism is given by the identification $\mu = [\xi_\bC] - [\underline{\bC}]$
\end{proposition}

Using the above identification of rings, we can show the following property relating the ring structure to the complex Radon-Hurwitz numbers:

\begin{proposition}\label{roc} For any positive integer $n\in\bN$, 
$$n\mu=0\mbox{ in }\widetilde{K}_\bC(\bR\bP^{d-1})\, \Longleftrightarrow \, d\leq \rho_\bC(n).$$
\end{proposition}
\begin{proof} First of all, remark that relation \ref{rel2} can be equivalently written (using relation \ref{rel1}) as $2^{g(d)}\mu=0$ so that $n\mu=0$ if and only if $n$ is an integer multiple of $2^{g(d)}$. Write now $n=2^a(2k+1)$ with $a,k\in\bN$, then
$$d\leq\rho_\bC(n)=2a+2\; \Leftrightarrow\; a\geq\frac{d-1}{2}-\frac{1}{2} \; \Leftrightarrow \; a \geq g(d),$$
where the last equivalence is given by the fact that $a$ is an integer and $g(d)$ is either $(d-1)/2$ or $(d-1)/2 -1/2$ according to the parity of $d$.
In conclusion,
$$n\mu=0 \; \Leftrightarrow\; 2^a(2k+1)\mu=0 \; \Leftrightarrow\; a\geq g(d) \; \Leftrightarrow\; d\leq\rho_\bC(n)$$
\end{proof}

\subsection{Odd maps} 

\begin{definition}
Any map $\psi:S^{d-1}\rightarrow \cM_{n}(\bC)$ verifying the relation $\psi(-x)=-\psi(x)$ is called an {\em odd map}.
\end{definition}

\begin{proposition}\label{oddmap}
Any odd map $\psi$ induces a morphism of vector bundles over $\bR\bP^{d-1}$:
$$\Psi:\underline{\bC}^n\rightarrow n\xi_\bC\quad\quad\mbox{ defined locally as }\Psi([x],v)=([x], \psi(x)v).$$
Moreover, if $\psi(x)$ has rank $r$ for any $x$, then $K=\ker \Psi$ and $C=\mbox{\em coker}\,\Psi$ are well-defined vector bundles
 and the following isomorhism holds:
 \begin{equation}\label{iso}
 K\oplus n\xi_\bC=\underline{\bC}^n\oplus C.
 \end{equation}
\end{proposition}
\begin{proof} Consider the map $$\Psi^\prime:S^{d-1}\times \bC^n\rightarrow S^{d-1}\times\bC^n\quad\quad\Psi^\prime(x,v)=(x,\psi(x)v).$$
Since $\psi$ is odd, $\Psi^\prime$ is equivariant with respect to the actions of $\bR^\ast$ on $S^{d-1}\times\bC^n$ given by 
$$f_\lambda(x,v)=\left(\frac{\lambda}{|\lambda |}x, v\right),\quad\quad g_\lambda(x,v)=\left(\frac{\lambda}{|\lambda |}x, \lambda v\right),\quad\lambda\in\bR^\ast.$$ The map induced by $\Psi^\prime$ by passing to the quotients is exactly $\Psi,$ indeed: $$\frac{S^{d-1}\times\bC^n}{f_\lambda}=\bR\bP^{d-1}\times\bC^n \quad\mbox{ and }\quad \frac{S^{d-1}\times\bC^n}{g_\lambda}=n(\xi\oplus\xi)=n\xi_\bC.$$
Isomorphism \ref{iso} is a consequence of the fact that any exact sequence of vector bundles splits.
\end{proof}

With this setting, we can prove the main theorem of Adams, Lax and Phillips in \cite{adams2}
\begin{theorem} 
If $V$ is a real vector space such that $V\setminus\{0\}\subset GL_n(\bC)$, then $\dim V\leq \rho_\bC(n).$
\end{theorem}
\begin{proof} Apply proposition \ref{oddmap} to the following setting: $d=\dim V$, $S^{d-1}$ the unit sphere of $V$ and $\psi$ the induced inclusion (wich trivially is an odd map) $$S^{d-1}\subset V\setminus\{0\}\subset GL_n(\bC)\subset\cM_n(\bC).$$
In this case $\psi(x)$ is always an invertible matrix, then isomorphism \ref{iso} becomes $n\xi_\bC=\underline{\bC}^n.$ In the ring $\widetilde{K}_\bC(\bR\bP^{d-1})$ this means that $n$ times the generator $\mu$ is zero, hence, by proposition \ref{roc}, $d\leq\rho_\bC(n).$
\end{proof}

\section{Adjugate matrices}
For a square matrix $A\in \cM_n(\bC)$, let $A_{i,j}$ be the submatrix obtained from $A$ by deleting its $i$-th row and $j$-th column and denote by $A^c$ the transpose of the adjugate of $A:$ 
$$(A^c)_{i,j}= (-1)^{i+j}\det A_{i,j}.$$

We define a map 
$$\psi:\cM_n(\bC)\rightarrow \cM_n(\bC)\quad \mbox{ by setting } \quad \psi(A)=A+i\overline{A^c},$$
the bar denoting complex conjugation. It is clear from the definition of $A^c$ that if $\mbox{rank} A\leq n-2$, then $\psi(A)=A$. Now set: 
$$\cZ=\{A\in\cM_n(\bC)| \mbox{ rank }A\geq n-1, \det A \neq ir, r<0 \}$$ 
and remark that both the spaces $H$ and $R$ (of hermitian and real $n\times n$ matrices) are subsets of $\cZ$.

\begin{proposition} For $\psi$ as above, $\psi(\cZ)\subset GL_n(\bC)$
\end{proposition}
\begin{proof} By contradiction, take $A\in\cZ$ and $v\in\bC^n, v\neq 0$ such that 

\begin{equation}\label{eq1} \psi(A)v= Av+i\overline{A^c}v=0.\end{equation}

\noindent Multiplying this equation on the left by $^t\overline{A}$ gives

\begin{equation}\label{eq2} ^t\overline{A}Av + i \det(\overline{A})v=0\end{equation}

\noindent Then, since $v\neq 0$, equation \ref{eq2} says that $-i\det(\overline A)$ is an eigenvector of the hermitian matrix $^t\overline{A}A$ and is {\em a fortiori } real and not negative. This means that $\det(A)=ir$ with $r\leq 0$ and, since $A\in \cZ$, the equality $\det(A)=0=\det(\overline{A})$ holds. Equation \ref{eq2} becomes $^t\overline{A}Av=0$ and, after left multiplication by $^t\overline{v}$,

$$ ^t\overline{v}\ ^t\overline{A}Av= \|Av\|=0.
$$

\noindent As a consequence, equation \ref{eq1} reads 

\begin{equation}\label{eq3} Av=\overline{A^c}v=0.
\end{equation}

\noindent From the fundamental property of the adjugate matrix $ A\ ^tA^c= \det(A) I =0$ it follows that the image of $^tA^c$ is contained in the null space of $A$ which is, by equation \ref{eq3}, spanned by $v$. Then, $\overline{A^c}v=0$ implies $$\overline{A^c}\ ^tA^c=0,$$
that is $A^c=0$ but this is impossible since $\mbox{rank}A\geq n-1$ forces at least one element of $A^c$ to be different from zero. 
\end{proof}

\begin{remark} Set $X=\{A\in\cM_n(\bC)| \mbox{ rank } A\leq n-2\}$. It is somehow interesting to observe that the family of maps $$\psi_s : \cM_n(\bC)\rightarrow \cM_n(\bC) \quad \quad \psi_s(A)=A+s\ i \overline{A^c}, \, s\in \bR $$ realizes, for $0\leq s\leq 1$ a homotopy from the identity $id(A)=A$ to $\psi=\psi_1$, fixing the algebraic variety $X$ pointwise. Moreover, for any $s>0$, $\psi_s (\cZ)\subset GL_n(\bC)$.
\end{remark}
\medskip

Multilinearity of the determinant imples that:
\begin{proposition} If $n$ is even, then $\psi(-A)=-\psi(A).$
\end{proposition}
\begin{proof} It is a direct consequence of the fact that the submatrices $A_{i,j}$ have odd order hence $\det(-A_{i,j})=-\det(A_{i,j})$ and consequently $(-A)^c=-A^c.$
\end{proof}

We can now prove:
\begin{proposition} If $V$ is a real vector space such that $V\setminus \{0\}\subset \cZ$ and $n$ is even, then 
$$\dim V \leq \rho_\bC(n).$$
\end{proposition}
\begin{proof} The restriction of $\psi$ to the unit sphere of $V$ is an odd map 
$$\psi:S^{d-1}\rightarrow GL_n(\bC)\quad \quad d=\dim V;$$ 
by proposition \ref{oddmap}, $\psi$ induces an isomorphism of vector bundles over $\bR\bP^{d-1}$: $$\underline{\bC}^n=n\xi_\bC.$$
In the ring of reduced complex K-theory $\widetilde{K}_\bC(\bR\bP^{d-1})$ this isomorphim means that $n$ times the generator $\mu$ is zero, hence proposition \ref{roc} implies $d\leq\rho_\bC(n).$
\end{proof}

As a corollary it follows:
\begin{theorem}\label{teo} If $n$ is even, then $$h_{n,n-1}\ \leq \ \rho_\bC(n) \quad\mbox{ and }\quad r_{n,n-1}\\ \leq \ \rho_\bC(n).$$
\end{theorem}
\begin{proof} This is due to the remark that both $H$ and $R$ are subsets of $\cZ.$
\end{proof}

\subsection{Some remarks on the estimates}
Theorem \ref{teo} provides a very sharp estimate on the numbers $h_{n,n-1}$ and $r_{n,n-1}$. This is particularly clear when we compare our result with the values of $h_{n,n}$ and $r_{n,n}$ given in \cite{adams2}. Indeed, for hermitian matrices, if we write as customary $n=2^a(2k+1)$ with $a,k\in\bN$, we get 
$$h_{n,n}=\rho_\bC(n/2)+1=2a+1\ \leq \ h_{n,n-1}\ \leq \ 2a+2=\rho_\bC(n) \quad \mbox{($n$ even)}.$$
For real matrices, on the other hand, we similarly have 
$$ r_{n,n}=\rho(n)\ \leq \ r_{n,n-1} \ \leq \ \rho_\bC(n)\quad \mbox{($n$ even)}.$$ It is interesting to compare  the values of $\rho(n)$ to those of $\rho_\bC(n)$: writing $n=2^{a+4b}(2k+1)$, with $a,b,k\in\bN$ and $0\leq a \leq 3$, we get the following
\begin{center}
\begin{tabular}{lcr}
$a$ & $\rho(n)$ & $\rho_\bC(n)$ \\
\hline 0 & 1+8b & 2+8b\\
1 & 2 +8b & 4 +8b \\
2 & 4+8b & 6+8b \\
3 & 8+8b & 8+8b
 \end{tabular}
\end{center}

thus, in particular, 

\begin{proposition}
If $8$ divides $n$ then $r_{n,n-1}=\rho(n)$.
\end{proposition}
This proposition gives a negative answer to a question posed in \cite{pacific} asking whether the inequalities in the sequence $$ r_{n,n}\leq r_{n,n-1}\leq\dots\leq r_{n,2}\leq r_{n,1} $$ were always sharp or not.

%%%%%%%%%%%%%%%%%%%%%%%%%%%%%%%%%%%%%%%%%%%%%%%%%%%%%%%%%%%%%%%%%%%%%%%%%%%%%%%%%%%%%%%%%%%%%%%%%
%%%%%%%%%%%%%%%%%%%%%%%%%%%%%%% BIBLIOGRAFIA %%%%%%%%%%%%%%%%%%%%%%%%%%%%%%%%%%%%%%%%%%%%%%%%%%%%
%%%%%%%%%%%%%%%%%%%%%%%%%%%%%%%%%%%%%%%%%%%%%%%%%%%%%%%%%%%%%%%%%%%%%%%%%%%%%%%%%%%%%%%%%%%%%%%%%

\vfill
\noindent
 {\sc Andrea Causin}\\
D.A.P., Universit\`a di Sassari\\
Piazza Duomo 6, 07041 Alghero (SS), Italia\\
{\tt acausin@uniss.it}\\

\end{document}